\documentclass[12pt,a4paper]{amsart}
\usepackage{amsmath, amsfonts,amssymb, amsthm, amscd}
\newcommand{\Z}{\mathbb Z}

\newtheorem{theo}{Theorem}

\title{On  presentations of generalizations of braids with few
generators}
\author{V.~Vershinin}
\address{D\'epartement des Sciences Math\'ematiques,
Universit\'e Montpellier II,
Place Eug\'ene Bataillon,
34095 Montpellier cedex 5, France}
\email{ vershini@math.univ-montp2.fr}
\address{ Sobolev Institute of Mathematics, Novosibirsk, 630090,
Russia }
\email{ versh@math.nsc.ru}
\keywords{Braid group, singular braid monoid, generalized braid groups,
presentation.}
\thanks{
The author was supported
 in part by the ACI project ACI-NIM-2004-243 "Braids and Knots". }
\begin{document}

\subjclass[2000]{Primary 20F36; Secondary 20F38, 20M05}

\begin{abstract} In his initial paper on braids E.~Artin gave
a presentation with two generators for an arbitrary braid group.
We give analogues of this Artin's presentation for various
generalizations of braids.
\end{abstract}
\maketitle

The diverse aspects of presentations of braid groups and their
generalizations continue to attract attention \cite{Ser}, \cite{BKL},
\cite{V}.

The canonical presentation of the braid group $Br_n$ was given by
E.~Artin \cite{Art1} and is well known. It has the generators
$\sigma_1$, $\sigma_2$, $\dots$, $\sigma_{n-1}$ and relations
\begin{equation*}
 \begin{cases} \sigma_i \sigma_j &=\sigma_j \, \sigma_i, \ \
\text{if} \ \ |i-j|>1, \ \ i,j =1, ..., n-1;
\\ \sigma_i \sigma_{i+1} \sigma_i &= \sigma_{i+1} \sigma_i \sigma_{i+1},
\ \ i =1, ..., n-2.
\end{cases}
\end{equation*}
There exist other presentations of the braid group.
J.~S.~Birman, K.~H.~Ko and S.~J.~Lee \cite{BKL} introduced the
presentation with generators $a_{ts}$ with $1 \leq s<t\leq n$
and relations
\begin{equation*}
\begin{cases}
a_{ts}a_{rq}&=a_{rq}a_{ts} \ \ {\rm for} \ \ (t-r)(t-q)(s-r)(s-q)>0,\\
a_{ts}a_{sr} &=a_{tr}a_{ts}=a_{sr}a_{tr}  \ \ {\rm for} \ \
1\leq r<s<t\leq n .
\end{cases}
\end{equation*}
The generators $a_{ts}$ are expressed in the canonical generators
$\sigma_i$ as follows:
 $$a_{ts}=(\sigma_{t-1}\sigma_{t-2}\cdots\sigma_{s+1})\sigma_s
(\sigma^{-1}_{s+1}\cdots\sigma^{-1}_{t-2}\sigma^{-1}_{t-1})  \ \
{\rm for} \ \ 1\leq s<t\leq n.$$
An analogue of the Birman-Ko-Lee presentation for the generalizations
of braids (namely, for the singular braid monoid) was obtained in \cite{V}.

In the initial paper \cite{Art1} Artin gave another presentation of
the braid group, with two generators, say $\sigma_1$ and $\sigma$,
and the following relations:
\begin{equation}
 \begin{cases}
\sigma_1 \sigma^i \sigma_1 \sigma^{-i} &=
\sigma^i \sigma_1 \sigma^{-i} \sigma_1 \ \  \text{for} \ \
2 \leq i\leq {n / 2}, \\
\sigma^n &= (\sigma \sigma_1)^{n-1}.
\end{cases} \label{eq:2relations}
\end{equation}

The connection with the canonical generators is given by the formulae:
 \begin{equation}
\sigma = \sigma_1 \sigma_{2} \dots \sigma_{n-1},
\label{eq:sigma}
\end{equation}
\begin{equation}
\sigma_{i+1} =\sigma^i \sigma_1 \sigma^{-i}, \quad i =1, \dots
{n-2}.
\label{eq:sigma_i}
\end{equation}

This presentation was also discussed in the book by H.~S.~M.~Coxeter
and W.~O.~J.~Moser \cite{CM}.

It is interesting to obtain the analogues of the presentations of the type
(\ref{eq:2relations}) for various generalizations of braids
\cite{Bri1}, \cite{Del}, \cite{Bae}, \cite{Bir2}, \cite{FRR2},
\cite{Ve4}.

Let us consider the braid group of a sphere $Br_n(S^2)$. It has the
presentation with generators $\delta_i$, $i=1, ..., n-1$ and relations:
\begin{equation*}
 \begin{cases} \delta_i \delta_j &=\delta_j \, \delta_i, \ \
\text{if} \ \ |i-j|
>1,\\
\delta_i \delta_{i+1} \delta_i &= \delta_{i+1} \delta_i \delta_{i+1}, \\
\delta_1 \delta_2 \dots \delta_{n-2}\delta_{n-1}^2\delta_{n-2} \dots
\delta_2\delta_1 &=1.
\end{cases}
\end{equation*}
This presentation was find by O.~Zariski \cite{Za1} in 1936 and
rediscovered later by E.~Fadell and J.~Van Buskirk \cite{FaV} in 1961.
From this presentation one can easily obtain the presentation with
two generators $\delta_1$, $\delta$ and the following relations:
\begin{equation*}
 \begin{cases}
\delta_1 \delta^i \delta_1 \delta^{-i} &=
\delta^i \delta_1 \delta^{-i} \delta_1 \ \  \text{for} \ \
2 \leq i\leq {n / 2}, \\
\delta^n &= (\delta \delta_1)^{n-1},\\
\delta^n  (\delta_1 \delta^{-1})^{n-1}&= 1.\\
\end{cases}
\end{equation*}

Another generalization of braids is the {\it Baez--Birman monoid}
$SB_n$ which is also called as {\it singular braid monoid}
\cite{Bae}, \cite{Bir2}. It is defined as a monoid with generators
$\sigma_i, \sigma_i^{-1}, x_i$, $i=1,\dots,n-1,$ and relations
\begin{equation}
\begin{cases}
&\sigma_i\sigma_j=\sigma_j\sigma_i, \ \text {if} \ \ |i-j| >1,\\
&x_ix_j=x_jx_i, \ \text {if} \ \ |i-j| >1,\\
&x_i\sigma_j=\sigma_j x_i, \
\text {if} \ \ |i-j| \not=1,\\
&\sigma_i \sigma_{i+1} \sigma_i = \sigma_{i+1} \sigma_i
\sigma_{i+1},\\
&\sigma_i \sigma_{i+1} x_i = x_{i+1} \sigma_i \sigma_{i+1},\\
&\sigma_{i+1} \sigma_ix_{i+1} = x_i \sigma_{i+1} \sigma_i,\\
&\sigma_i\sigma_i^{-1}=\sigma_i^{-1}\sigma_i =1.
\end{cases}
\label{eq:singrel}
\end{equation}
In pictures, $\sigma_i$ corresponds to the canonical generator of
the braid group and $x_i$ represents an intersection
of the $i$th and $(i+1)$st strand as in
Figure~\ref{fi:singene}. The singular braid monoid on two strings
$SB_2$ is isomorphic to $\Z\oplus\Z^+$.

\begin{figure}
\begin{picture}(0,130)(0,-10) 
\thicklines
\put(0,50){\circle*{5}} \put(-100,100){\line(0,-1){100}}
\put(-50,100){\line(0,-1){100}} \put(-25,100){\line(1,-2){50}}
\put(25,100){\line(-1,-2){50}} \put(50,100){\line(0,-1){100}}
\put(100,100){\line(0,-1){100}}
\put(-100,110){\makebox(0,0)[cc]{$1$}}
\put(-50,110){\makebox(0,0)[cc]{$i-1$}}
\put(-25,110){\makebox(0,0)[cc]{$i$}}
\put(25,110){\makebox(0,0)[cc]{$i+1$}}
\put(50,110){\makebox(0,0)[cc]{$i+2$}}
\put(100,110){\makebox(0,0)[cc]{$n$}}
\put(-75,50){\makebox(0,0)[cc]{.\quad.\quad.}}
\put(75,50){\makebox(0,0)[cc]{.\quad.\quad.}}
\end{picture}

\caption{}\label{fi:singene}
\end{figure}

Motivation for the introduction of this object was the
Vassiliev -- Goussarov theory of finite type invariants.

Let $F_n$ be a free group on $n$ generators $x_1, \dots, x_n$ and
$\operatorname{Aut} F_n$
its automorphism group. The {\it braid-permutation group} $BP_n$,
considered by R.~Fenn,
R.~Rim\'anyi and C.~Rourke \cite{FRR2}, is
the subgroup of $\operatorname{Aut} F_n$, generated
by two sets of the automorphisms: $\sigma_i$
\begin{equation*} \begin{cases}
x_i &\mapsto x_{i+1},
\\ x_{i+1} &\mapsto x_{i+1}^{-1}x_ix_{i+1}, \\
x_j &\mapsto x_j, j\not=i,i+1,
\end{cases}
\end{equation*}
and  $\xi_i$:
\begin{equation*}
\begin{cases} x_i &\mapsto x_{i+1}, \\ x_{i+1} &\mapsto
x_i,     \\ x_j &\mapsto x_j, \ \ j\not=i,i+1. \end{cases}
\end{equation*}
R.~Fenn,
R.~Rim\'anyi and C.~Rourke proved \cite{FRR2} that
this group
is given by the set of generators: $\{ \xi_i, \sigma_i, \ \ i=1,2,
..., n-1 \}$ and relations:
$$ \begin{cases} \xi_i^2&=1, \\ \xi_i \xi_j
&=\xi_j \xi_i, \ \ \text {if} \ \ |i-j| >1,\\ \xi_i \xi_{i+1}
\xi_i &= \xi_{i+1} \xi_i \xi_{i+1}. \end{cases} $$
\vglue 0.01cm
\centerline { The symmetric group relations}
$$ \begin{cases} \sigma_i \sigma_j &=\sigma_j
\sigma_i, \ \text {if} \  |i-j| >1,
\\ \sigma_i \sigma_{i+1} \sigma_i &= \sigma_{i+1} \sigma_i \sigma_{i+1}.
\end{cases} $$
\vglue 0.01cm
\centerline {The braid group relations } $$ \begin{cases} \sigma_i \xi_j
&=\xi_j \sigma_i, \ \text {if} \  |i-j| >1,
\\ \xi_i \xi_{i+1} \sigma_i &= \sigma_{i+1} \xi_i \xi_{i+1},
\\ \sigma_i \sigma_{i+1} \xi_i &= \xi_{i+1} \sigma_i \sigma_{i+1}.
\end{cases} $$
\vglue0.01cm
\centerline {The mixed relations }
\smallskip

R.~Fenn, R.~Rim\'anyi and C.~Rourke also gave a geometric
interpretation of $BP_n$ as a group of welded braids.

Let us obtain a presentation of the singular braid monoid and the
braid-permutation group analogous to (\ref{eq:2relations}).
If we add the new generator $\sigma$,
defined by (\ref{eq:sigma}) to the set of generators of $SB_n$ then the following
relations hold
\begin{equation}
x_{i+1} =\sigma^i x_1 \sigma^{-i}, \quad i =1, \dots
{n-2}.
\label{eq:x_i}
\end{equation}
This gives a possibility to get rid of $x_i$, $i\geq 2$.

\begin{theo} The singular braid monoid $SB_n$ has a presentation with
generators $\sigma_1$, $\sigma_1^{-1}$, $\sigma$, $\sigma^{-1}$ and
$x_1$
and relations
\begin{equation}
 \begin{cases}
\sigma_1 \sigma^i \sigma_1 \sigma^{-i} =
\sigma^i \sigma_1 \sigma^{-i} \sigma_1 \ \  \text{for} \ \
2 \leq i\leq {n / 2}, \\
\sigma^n = (\sigma \sigma_1)^{n-1},\\
x_1\sigma^i\sigma_1\sigma^{-i}= \sigma^i\sigma_1\sigma^{-i} x_1
\ \  \text{for} \ \ i=0, 2, \dots, {n - 2}, \\
x_1 \sigma^i x_1 \sigma^{-i} =
\sigma^i x_1 \sigma^{-i} x_1 \ \  \text{for} \ \
2 \leq i\leq {n / 2}, \\
\sigma^n x_1 = x_1\sigma^n,\\
x_1\sigma\sigma_1\sigma^{-1}\sigma_1 = \sigma\sigma_1\sigma^{-1}\sigma_1
\sigma x_1\sigma^{-1}, \\
\sigma_1\sigma_1^{-1}=\sigma_1^{-1}\sigma_1 =1,\\
\sigma\sigma^{-1}=\sigma^{-1}\sigma =1.
\end{cases} \label{eq:2gens}
\end{equation}
\end{theo}
\begin{proof} We follow the original Artin's proof \cite{Art1} and
we begin with the presentation of $SB_n$ using the generators
$\sigma_i, \sigma_i^{-1}, x_i$, $i=1,\dots,n-1,$ and relations
(\ref{eq:singrel}). Then we add the new generators
$\sigma$, $\sigma^{-1}$,  relation (\ref{eq:sigma}) and the following
relations
\begin{equation*}
\sigma\sigma^{-1}=\sigma^{-1}\sigma =1.
\end{equation*}
Consider $\sigma x_i$. Using the braid relations in the same way as Artin
considered $\sigma\sigma_i$ we have
\begin{multline*}
\sigma x_i= \sigma_1 \dots \sigma_{n-1} x_i =
\sigma_1 \dots\sigma_i\sigma_{i+1}x_i
\sigma_{i+2}\dots\sigma_{n-1} = \\
=\sigma_1 \dots x_{i+1}\sigma_i\sigma_{i+1}
\sigma_{i+2}\dots\sigma_{n-1} = x_{i+1}\sigma.
\end{multline*}
We arrive thus at relations (\ref{eq:sigma_i}) and (\ref{eq:x_i}).
Now we can get rid from the fifth relation in
(\ref{eq:singrel}). First of all using (\ref{eq:sigma_i}) and (\ref{eq:x_i})
it is reduced to the case $i=1$,
which is considered as follows: $x_1$ commutes with
$\sigma_2^{-1}\sigma_1^{-1} \sigma$:
\begin{equation*}
x_1 \sigma_2^{-1}\sigma_1^{-1} \sigma =
\sigma_2^{-1}\sigma_1^{-1} \sigma x_1 =
\sigma_2^{-1}\sigma_1^{-1}x_2 \sigma .
\end{equation*}
So, we obtain
\begin{equation*}
x_1 \sigma_2^{-1}\sigma_1^{-1}  = \sigma_2^{-1}\sigma_1^{-1}x_2.
\end{equation*}
This is equivalent to the fifth relation in (\ref{eq:singrel}) for $i=1$.
Using (\ref{eq:sigma_i}) and (\ref{eq:x_i}) the sixth relation in
(\ref{eq:singrel}) is reduced to the case $i=1$, which is the sixth relation
in (\ref{eq:2gens}).

The third and forth relations in (\ref{eq:2gens}) are easy consequences
of the corresponding relations in (\ref{eq:singrel}) and (\ref{eq:sigma_i})
and (\ref{eq:x_i}). The fifth relation in (\ref{eq:2gens}) is a
consequence of the definition of $\sigma$ and the singular braid relations
(\ref{eq:singrel}).
The third relation in (\ref{eq:singrel}) is obtained from the forth and
fifth relations in (\ref{eq:2gens}) in the same way as Artin obtained the
commutation
of $\sigma_i$ and $\sigma_j$ from the relations in (\ref{eq:2relations}).
Essentially, Artin used the fact that it follows from the second relation
in (\ref{eq:2relations}) that $\sigma^n$ is in the center of $Br_n$. Here
we need the fifth relation in (\ref{eq:2gens}) to have this fact in $SB_n$.
The third relation in (\ref{eq:singrel}) in the new generators is rewritten
as follows
\begin{equation*}
\sigma^ix_1\sigma^{-i} \sigma^j \sigma_1 \sigma^{-j} =
\sigma^j \sigma_1 \sigma^{-j}\sigma^i x_1\sigma^{-i},
\end{equation*}
what is equivalent to
\begin{equation*}
 x_1\sigma^{j-i}  \sigma_1 \sigma^{i-j} =
\sigma^{j-i} \sigma_1 \sigma^{i-j} x_1.
\end{equation*}
If $j>i$ then this is exactly the third relation in (\ref{eq:2gens}), if
$j<i$ then it follows from the third relation in (\ref{eq:2gens}) by
conjugation
by $\sigma^n$ and using the commutation of $\sigma^n$ with $x_1$.
\end{proof}

For the case of the braid-permutation group $SB_n$ we also add the new
generator
$\sigma$, defined by (\ref{eq:sigma}) to the set of standard generators of
$BP_n$; then relations (\ref{eq:sigma_i}) and the following relations hold
\begin{equation*}
\xi_{i+1} =\sigma^i \xi_1 \sigma^{-i}, \quad i =1, \dots,
{n-2}.
\end{equation*}
This gives a possibility to get rid of $\xi_i$ as well as of $\sigma_i$
for $i\geq 2$.

\begin{theo} The braid-permutation group $BP_n$ has presentation with
generators $\sigma_1$, $\sigma$,  and $\xi_1$ and relations
\begin{equation*}
 \begin{cases}
\sigma_1 \sigma^i \sigma_1 \sigma^{-i} =
\sigma^i \sigma_1 \sigma^{-i} \sigma_1 \ \  \text{for} \ \
2 \leq i\leq {n / 2}, \\
\sigma^n = (\sigma \sigma_1)^{n-1},\\
\xi_1\sigma^i\sigma_1\sigma^{-i}= \sigma^i\sigma_1\sigma^{-i} \xi_1
\ \  \text{for} \ \ i= 2 \dots {n - 2}, \\
\xi_1 \sigma^i \xi_1 \sigma^{-i} =
\sigma^i \xi_1 \sigma^{-i} \xi_1 \ \ \text{for} \ \ i= 2 \dots {n - 2}, \\
\xi_1\sigma\xi_1\sigma^{-1}\sigma_1 = \sigma\sigma_1\sigma^{-1}\xi_1
\sigma \xi_1\sigma^{-1}, \\
\xi_1\sigma\xi_1\sigma^{-1}\xi_1 = \sigma\xi_1\sigma^{-1}\xi_1
\sigma \xi_1\sigma^{-1}, \\
\xi^2=1.
\end{cases}
\end{equation*}     $\square$
\end{theo}

Let us consider the generalized braid groups in the sense of Brieskorn
(or so the called Artin groups) \cite{Bri1}.
It is easy to see that for the braid groups of type $B_n$
from the canonical presentation
with generators  $\sigma_i$, $i=1, \dots, n-1$ and $\tau$, and relations:
\begin{equation*}
\begin{cases}
\sigma_i\sigma_j &=\sigma_j\sigma_i, \ \
\text{if} \ \ |i-j| >1,
\\ \sigma_i \sigma_{i+1} \sigma_i &=
\sigma_{i+1} \sigma_i \sigma_{i+1},\\
\tau\sigma_i &=\sigma_i\tau, \ \ \text{if} \ \  i\geq 2,\\
\tau\sigma_1\tau\sigma_1&=\sigma_1\tau\sigma_1\tau , \\
\end{cases}
\end{equation*}
we can obtain the presentation with three generators
$\sigma_1$, $\sigma$ and $\tau$ and the following relations:
\begin{equation}
 \begin{cases}
\sigma_1 \sigma^i \sigma_1 \sigma^{-i} &=
\sigma^i \sigma_1 \sigma^{-i} \sigma_1 \ \  \text{for} \ \
2 \leq i\leq {n / 2}, \\
\sigma^n &= (\sigma \sigma_1)^{n-1},\\
\tau\sigma^i\sigma_1\sigma^{-i} &=\sigma^i\sigma_1\sigma^{-i}\tau
\ \ \text{for} \ \ 2 \leq i\leq {n - 2}, \\
\tau\sigma_1\tau\sigma_1&=\sigma_1\tau\sigma_1\tau . \\
\end{cases} \label{eq:2relB}
\end{equation}

If we add the following relations
\begin{equation*}
 \begin{cases}
\sigma_1^2 &= 1, \\
\tau^2 &= 1  \\
\end{cases}
\end{equation*}
to (\ref{eq:2relB}) we then arrive at a presentation of the Coxeter
group of type $B_n$.

Similarly, for the braid groups of the type $D_n$
from the canonical presentation
with generators  $\sigma_i$ and $\rho$, and relations:
\begin{equation*}
\begin{cases}
\sigma_i\sigma_j &=\sigma_j\sigma_i \ \
\text{if} \ \ |i-j| >1,\\ \sigma_i \sigma_{i+1} \sigma_i &=
\sigma_{i+1} \sigma_i \sigma_{i+1},\\
\rho\sigma_i &=\sigma_i\rho \ \ \text{if} \ \  i =1, 3, \dots, n-1,\\
\rho\sigma_2\rho&=\sigma_2\rho\sigma_2 , \\
\end{cases}
\end{equation*}
we can obtain the presentation with three generators
$\sigma_1$, $\sigma$ and $\rho$ and the following relations:
\begin{equation}
 \begin{cases}
\sigma_1 \sigma^i \sigma_1 \sigma^{-i} &=
\sigma^i \sigma_1 \sigma^{-i} \sigma_1 \ \  \text{for} \ \
2 \leq i\leq {n / 2}, \\
\sigma^n &= (\sigma \sigma_1)^{n-1},\\
\rho\sigma^i\sigma_1\sigma^{-i} &=\sigma^i\sigma_1\sigma^{-i}\rho
\ \ \text{for} \ \ i=0, 2, \dots,  {n - 2}, \\
\rho\sigma\sigma_1\sigma^{-1}\rho&=\sigma\sigma_1\sigma^{-1}
\rho\sigma\sigma_1\sigma^{-1} . \\
\end{cases} \label{eq:2relD}
\end{equation}

If we add the following relations
\begin{equation*}
 \begin{cases}
\sigma_1^2 &= 1, \\
\rho^2 &= 1  \\
\end{cases}
\end{equation*}
to (\ref{eq:2relD}) we come to a presentation of the Coxeter
group of type $D_n$.

For the exceptional braid groups of types $E_6 - E_8$ our presentations
look similar to the presentation for the groups of type $D$
(\ref{eq:2relD}).
We give it here for $E_8$: it has three generators
$\sigma_1$, $\sigma$ and $\omega$ and the following relations:
\begin{equation}
 \begin{cases}
\sigma_1 \sigma^i \sigma_1 \sigma^{-i} &=
\sigma^i \sigma_1 \sigma^{-i} \sigma_1 \ \  \text{for} \ \
i = 2, 3 , 4, \\
\sigma^8 &= (\sigma \sigma_1)^{7},\\
\omega\sigma^i\sigma_1\sigma^{-i} &=\sigma^i\sigma_1\sigma^{-i}\omega
\ \ \text{for} \ \ i= 0, 1, 3, 4, 5, 6, \\
\omega\sigma^2\sigma_1\sigma^{-2}\omega &=
\sigma^2\sigma_1\sigma^{-2}\omega\sigma^2\sigma_1\sigma^{-2}.
\end{cases} \label{eq:2relE}
\end{equation}
Similarly, if we add the following relations
\begin{equation*}
 \begin{cases}
\sigma_1^2 &= 1, \\
\omega^2 &= 1  \\
\end{cases}
\end{equation*}
to (\ref{eq:2relE}) we arrive at a presentation of the Coxeter
group of type $E_8$.

As for the other exceptional braid groups, $F_4$ has four generators and
it follows from its Coxeter diagram that there is no sense to speak about
analogues of the Artin presentation (\ref{eq:2relations}), $G_2$ and
$I_2(p)$ already have two generators and $H_3$ has three generators. For
$H_4$
it is possible to diminish the number of generators from four to three and
the presentation will be similar to that of $B_4$.

We can summarize informally what we were doing. Let a group have
a presentation
which can be expressed by a ``Coxeter-like" graph. If there exists
a linear subgraph corresponding to the standard presentation of the
classical braid group, then in the ``braid-like" presentation of
our group the part that corresponds to the linear subgraph can be
replaced by two generators and relations (\ref{eq:2relations}).
This reáè§ã can be applied to the complex reflection groups
\cite{ShT} whose  ``Coxeter-like" presentations is obtained in
\cite{BMR}, \cite{BM}. For the series of the complex braid groups
$B(2e, e,r)$, $e\geq 2$, $r\geq 2$ which correspond to the complex
reflection groups $G(de, e, r)$, $d\geq 2$ \cite{BMR} we take the linear
subgraph with nodes $\tau_2, \dots, \tau_r$, and put as above
$\tau = \tau_2 \dots \tau_r$. The group $B(2e, e,r)$ have
presentation with generators $\tau_2, \tau$,
$\sigma, \tau_2^\prime$ and relations
\begin{equation}
 \begin{cases}
\tau_2 \tau^i \tau_2 \tau^{-i} &=
\tau^i \tau_2 \tau^{-i} \tau_2 \ \  \text{for} \ \
2 \leq i\leq {r / 2}, \\
\tau^r &= (\tau \tau_2)^{r-1}, \\
\sigma \tau^i\tau_2\tau^{-i} &= \tau^i\tau_2\tau^{-i}  \sigma,
\ \  \text{for} \ \ 1 \leq i\leq {r - 2}, \\
\sigma \tau_2^\prime\tau_2 &= \tau_2^\prime \tau_2 \sigma, \\
\tau_2^\prime \tau\tau_2\tau^{-1}\tau_2^\prime &=
\tau\tau_2\tau^{-1} \tau_2^\prime\tau\tau_2\tau^{-1}, \\
\tau\tau_2\tau^{-1}\tau_2^\prime \tau_2 \tau\tau_2\tau^{-1}\tau_2^\prime
\tau_2 &= \tau_2^\prime\tau_2 \tau\tau_2\tau^{-1} \tau_2^\prime \tau_2
 \tau\tau_2\tau^{-1}, \\
\underbrace{\tau_2 \sigma \tau_2^\prime\tau_2 \tau_2^\prime
\tau_2 \tau_2^\prime \dots}_{ e+1 \  \text{factors}}&=
 \underbrace{\sigma\tau_2^\prime  \tau_2
\tau_2^\prime \tau_2
\tau_2^\prime \tau_2 \dots}_{ e+1 \  \text{factors}} \ .
\end{cases} \label{eq:2relBde}
\end{equation}
If we add the following relations
\begin{equation*}
 \begin{cases}
\sigma^d &= 1, \\
\tau_2^2 &= 1,  \\
{{\tau_2}^\prime}^2 &= 1  \\
\end{cases}
\end{equation*}
to (\ref{eq:2relBde}) we come to a presentation of the complex
reflection group $G(de, e, r)$.

The braid group $B(d,1, n)$, $d>1$, has the same presentation as the
Artin -- Brieskorn group of type $B_n$, but if we
add the following relations
\begin{equation*}
 \begin{cases}
\sigma_1^2 &= 1, \\
\tau^d &= 1  \\
\end{cases}
\end{equation*}
to (\ref{eq:2relB}) then we arrive at a presentation of the complex
reflection group $G(d, 1, n)$, $d\geq 2$.

For the series of braid groups $B(e, e,r)$, $e\geq 2$, $r\geq 3$
which correspond to the complex
reflection groups $G(e, e, r)$, $e\geq 2$, $r\geq 3$
  we take again the linear
subgraph with the nodes $\tau_2, \dots, \tau_r$, and put as above
$\tau = \tau_2 \dots \tau_r$. The group $B(e, e,r)$ may have
the presentation with generators $\tau_2, \tau$,
$\tau_2^\prime$ and relations
\begin{equation}
 \begin{cases}
\tau_2 \tau^i \tau_2 \tau^{-i} &=
\tau^i \tau_2 \tau^{-i} \tau_2 \ \  \text{for} \ \
2 \leq i\leq {r / 2}, \\
\tau^r &= (\tau \tau_2)^{r-1}, \\
\tau_2^\prime \tau\tau_2\tau^{-1}\tau_2^\prime &=
\tau\tau_2\tau^{-1} \tau_2^\prime\tau\tau_2\tau^{-1}, \\
\tau\tau_2\tau^{-1}\tau_2^\prime \tau_2 \tau\tau_2\tau^{-1}\tau_2^\prime
\tau_2 &= \tau_2^\prime\tau_2 \tau\tau_2\tau^{-1} \tau_2^\prime \tau_2
 \tau\tau_2\tau^{-1}, \\
\underbrace{\tau_2 \tau_2^\prime\tau_2 \tau_2^\prime
\tau_2 \tau_2^\prime \dots}_{ e \  \text{factors}}&=
 \underbrace{\tau_2^\prime  \tau_2
\tau_2^\prime \tau_2
\tau_2^\prime \tau_2 \dots}_{ e \  \text{factors}}.
\end{cases} \label{eq:2relBe}
\end{equation}
If $e=2$ then this presentation is the same as for the presentaion
for the Artin -- Brieskorn group of type $D_r$ (\ref{eq:2relD}).
If we add the following relations
\begin{equation*}
 \begin{cases}
\tau_2^2 &= 1,  \\
{{\tau_2}^\prime}^2 &= 1  \\
\end{cases}
\end{equation*}
to (\ref{eq:2relBe}), then we obtain a presentation of the complex
reflection group $G(e, e, r)$, $e\geq 2$, $r\geq 3$.

As for the exceptional (complex) braid groups, it is reasonable to
consider the groups $Br(G_{30})$, $Br(G_{33})$ and
 $Br(G_{34})$  which correspond to the complex reflection groups
$G_{30}$, $G_{33}$ and $G_{34}$.

The presentation for $Br(G_{30})$ is similar to the presentation
(\ref{eq:2relB}) of $Br(B_4)$ with the last relation replaced by the
relation of length 5: the three generators
$\sigma_1$, $\sigma$ and $\tau$ and the following relations:
\begin{equation}
 \begin{cases}
\sigma_1 \sigma^2 \sigma_1 \sigma^{-2} &=
\sigma^2 \sigma_1 \sigma^{-2} \sigma_1 , \\
\sigma^4 &= (\sigma \sigma_1)^{3},\\
\tau\sigma^i\sigma_1\sigma^{-i} &=\sigma^i\sigma_1\sigma^{-i}\tau
\ \ \text{for} \ \ i = 2, 3, \\
\tau\sigma_1\tau\sigma_1\tau&=\sigma_1\tau\sigma_1\tau \sigma_1. \\
\end{cases} \label{eq:2relBG30}
\end{equation}
If we add the following relations
\begin{equation*}
 \begin{cases}
\sigma_1^2 &= 1, \\
\tau^2 &= 1  \\
\end{cases}
\end{equation*}
to (\ref{eq:2relBG30}), then we obtain a presentation of complex
reflection group $G_{30}$.

As for the groups $Br(G_{33})$ and $Br(G_{34})$,
we give here the presentation for
the letter one because the ``Coxeter-like" graph for $Br(G_{33})$
has one node less in the linear subgraph (discussed earlier) than that
of $Br(G_{34})$.
This presentation has the three generators
$s$, $z $ ($z=stuvx$ in the reflection generators) and $w$ and the
following relations:
\begin{equation}
 \begin{cases}
s z^i s z^{-i} &=
z^i s z^{-i} s \ \  \text{for} \ \
i = 2, 3 , \\
z^6 &= (z s)^{5},\\
w z^i s w^{-i} &=z^i s z^{-i}w
\ \ \text{for} \ \ i= 0, 3, 4, \\
w z^i s z^{-i} w &=
z^i s z^{-i} w z^i s z^{-i} \ \  \text{for} \ \
i = 1, 2, \\
w z^2 s z^{-2} w z s z^{-1} w z^2 s z^{-2}&=
 z s z^{-1} w z^2 s z^{-2} w z s z^{-1} w . \\
\end{cases} \label{eq:2relBG34}
\end{equation}
The same way if we add the following relations
\begin{equation*}
 \begin{cases}
s^2 &= 1, \\
w^2 &= 1  \\
\end{cases}
\end{equation*}
to (\ref{eq:2relBG34}), then we come to a presentation of the complex
reflection group $G_{34}$.

We can obtain presentations with few generators for
the other complex reflection groups using the already observed presentations
of the braid groups. For $G_{25}$ and $G_{32}$ we can use the
presentations (\ref{eq:2relations}) for the classical braid groups $Br_4$
and $Br_5$ with the only one additional relation
\begin{equation*}
 \begin{cases}
\sigma_1^3 &= 1. \\
\end{cases}
\end{equation*}

The author expresses his gratitude to V.~P.~Lexine for useful advices.

\end{document}